# Zero-Sum Triangles for Involutory, Idempotent, Nilpotent and Unipotent Matrices

Pengwei Hao[1], Chao Zhang[2], Huahan Hao[3]

**Abstract**: In some matrix formations, factorizations and transformations, we need special matrices with some properties and we wish that such matrices should be easily and simply generated and of integers. In this paper, we propose a zero-sum rule for the recurrence relations to construct integer triangles as triangular matrices with involutory, idempotent, nilpotent and unipotent properties, especially nilpotent and unipotent matrices of index 2. With the zero-sum rule we also give the conditions for the special matrices and the generic methods for the generation of those special matrices. The generated integer triangles are mostly newly discovered, and more combinatorial identities can be found with them.

## 1. Introduction

Zero-sum was coined in early 1940s in the field of game theory and economic theory, and the term "zero-sum game" is used to describe a situation where the sum of gains made by one person or group is lost in equal amounts by another person or group, so that the net outcome is neutral, or the sum of losses and gains is zero [1]. It was later also frequently used in psychology and political science. In this work, we apply the zero-sum rule to three adjacent cells to construct integer triangles.

Pascal's triangle is named after the French mathematician Blaise Pascal, but the history can be traced back in almost 2,000 years and is referred to as the Staircase of Mount Meru in India, the Khayyam triangle in Persia (Iran), Yang Hui's triangle in China, Apianus's Triangle in Germany, and Tartaglia's triangle in Italy [2]. In mathematics, Pascal's triangle is a triangular array of the binomial coefficients, symmetrically with both boundary edges being all 1s and every interior number being the sum of the two numbers directly above it. More generic $(a_n, b_n)$-Pascal's triangle [3] keeps the recurrence rule of the original Pascal's triangle unchanged for the interior cells of the triangle but uses integer sequences $a_n$ and $b_n$ for the input of the two boundary edges of the triangle.

With Pascal's triangle, there have been quite a few topics in mathematics related, a staggering number of combinatorial identities found [4], and a growing number of applications developed. Placing a Pascal's triangle into a matrix, we can use a unit lower triangular matrix to represent the triangle. In [5] for reversible resampling, we have found the matrix relations between Pascal's triangle matrix, Vandermonde matrix, matrices of the falling and rising factorial polynomials, and matrices of Stirling numbers of the first kind and the second kind.

Idempotent matrices, also known as projection matrices and closely related to involutory matrices, are used extensively in linear regression analysis [6] and econometrics [7], and involutory matrices, are used in quantum mechanics for Pauli spin matrices [8], and Hill cipher to avoid inversion, where the decryption needs the inverse of the key-matrix [9]. Idempotent matrices are diagonalizable and the eigenvalues are either 0 or 1 [10]. Involutory matrices are also diagonalizable and the eigenvalues are either $-1$ or 1 [11]. The matrices have many interesting properties [11]: If $\boldsymbol{T}$ is idempotent, $\boldsymbol{T}^2 = \boldsymbol{T}$, then $2\boldsymbol{T} - \boldsymbol{I}$ and $\boldsymbol{I} - 2\boldsymbol{T}$ are involutory. If $\boldsymbol{T}$ is involutory, $\boldsymbol{T}^2 = \boldsymbol{I}$, then we have: $\boldsymbol{T}^T\boldsymbol{T}$ and $\boldsymbol{T}\boldsymbol{T}^T$ are symmetric and orthogonal; $(\boldsymbol{I} + \boldsymbol{T})(\boldsymbol{I} - \boldsymbol{T}) = 0$; and $(\boldsymbol{I} + \boldsymbol{T})/2$ and $(\boldsymbol{I} - \boldsymbol{T})/2$ are idempotent. A matrix is similar to its inverse if and only if the matrix is the product of two involutory matrices, and a matrix of determinant 1 is the product of up to 4 involutory matrices.

Unipotent matrices, closely related to nilpotent matrices, are popularly used in linear transformations and matrix factorizations, e.g. unipotent lower triangular matrices for Gauss elimination and triangular factorizations [12]. Unipotent matrices and nilpotent matrices are not diagonalizable, and the eigenvalues of unipotent and nilpotent matrices are all 1s and all 0s, respectively [11]. A unipotent matrix $\boldsymbol{S}$ can be represented with a nilpotent matrix $\boldsymbol{N}$, $\boldsymbol{S} = \boldsymbol{I} + \boldsymbol{N}$, and we have $(\boldsymbol{S} - \boldsymbol{I})^k = \boldsymbol{N}^k = 0$, where $k$ is the index of the nilpotent matrix and the index of the unipotent matrix. For a unipotent matrix of index 2 and a nilpotent matrix of index 2, $k = 2$, we have $(\boldsymbol{S} - \boldsymbol{I})^2 = \boldsymbol{N}^2 = 0$ and $\boldsymbol{S}^{-1} = (\boldsymbol{I} + \boldsymbol{N})^{-1} = \boldsymbol{I} - \boldsymbol{N}$.

[1] Queen Mary University of London, E1 4NS, UK, p.hao@qmul.ac.uk
[2] Peking University, Beijing, 100871, China, chzhang@cis.pku.edu.cn
[3] University of Cambridge, Cambridge, CB2 1TL, UK, vincent.hao63@gmail.com

Fong and Sourour showed that a matrix of determinant 1 is the product of 3 unipotent matrices [13], and Wang and Wu proved that any unipotent matrix is the product of 2 unipotent matrices of index 2 and any matrix of determinant 1 is the product of 4 unipotent matrices of index 2 [14] or 4 unit triangular matrices (PLUS factorization) [15].

In this paper, we propose the zero-sum recurrence rule for triangle construction, take the triangles as triangular matrices to show the relations and properties, and investigate how the triangles as triangular matrices can be idempotent, involutory, nilpotent and unipotent.

We first introduce the zero-sum rule and the general term formula in Section 2, give two lemmas as the recursive formula for the squared triangle in Section 3 and the product of the triangle with the boundary in Section 4, in Section 5-8, we present four theorems for the methods to construct the triangles with various properties, and finally in Section 9 we wish to express some of our thoughts about the zero-sum triangles beyond this paper.

## 2. The Zero-Sum Rule and the Triangles

The triangular zero-sum rule, we propose, is to make zero the sum of three adjacent numbers for all the interior cells, every interior number and the two numbers directly above it, as shown in Figure 1, $t(i,k-1)+t(i,k)+t(i+1,k)=0$, where $i \geq k \geq 1$.

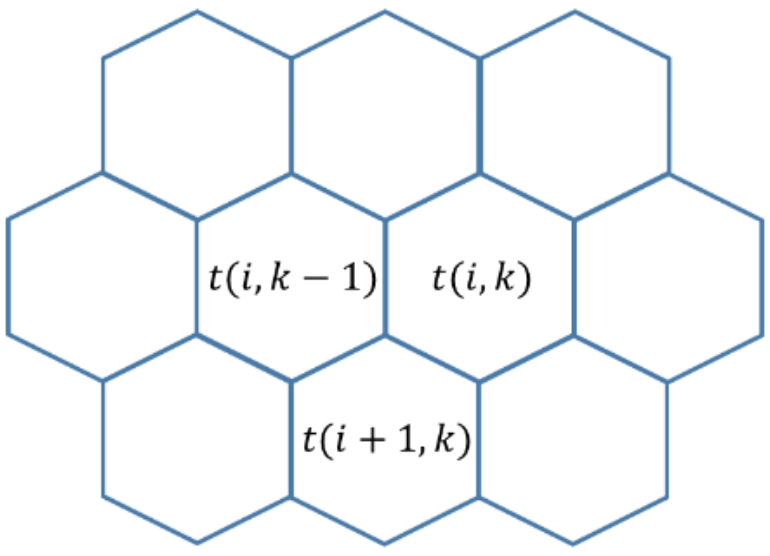

Figure 1. Three cells involved in zero-sum rule

See the two following triangles (1&2) for example of the zero-sum triangles:

```
                1
              1   0
            1  -1   1
          1   0   0   0
        1  -1   0   0   1
      1   0   1   0  -1   0
    1  -1  -1  -1   1   1   1
  1   0   2   2   0  -2  -2   0
1  -1  -2  -4  -2   2   4   2   1
```

Triangle 1

```
                1
              2  -1
            1  -1   1
          2   0   0  -1
        1  -2   0   1   1
      2   1   2  -1  -2  -1
    1  -3  -3  -1   3   3   1
  2   2   6   4  -2  -6  -4  -1
1  -4  -8 -10  -2   8  10   5   1
```

Triangle 2

If we use $t(i,k)$ to denote the $k$-th cell of the $i$-th row in the triangle, $i \geq k \geq 0$, the integer sequences of the two boundary edges can be given as $a_n$ and $b_n$, or:

$$\begin{gathered} t(0,0)=a_0=b_0, \\ t(1,0)=a_1, t(2,0)=a_2, \ldots, t(i,0)=a_i, \ldots, \\ t(1,1)=b_1, t(2,2)=b_2, \ldots, t(i,i)=b_i, \ldots \end{gathered} \tag{1}$$

To avoid involving exterior cells, the zero-sum rule can be equivalently expressed for three $(i,k)$ adjacent cells of the triangle as :

$$\begin{aligned} t(i-1,k-1)+t(i-1,k)+t(i,k)&=0 \qquad (i-1 \geq k \geq 1) \\ t(i,k-1)+t(i,k)+t(i+1,k)&=0 \qquad (i \geq k \geq 1) \\ t(i,k)+t(i,k+1)+t(i+1,k+1)&=0 \qquad (i-1 \geq k \geq 0) \end{aligned} \tag{2}$$

To find the $(i+1)$-th row with the $i$-th row, the recurrence rule for the interior cells of the triangle can also be written as:

$$t(i+1,k)=-t(i,k-1)-t(i,k) \qquad (i \geq k \geq 1) \tag{3}$$

Placing the $n$-row triangle into a lower triangular matrix, we have the $n$-row matrix as:

$$\boldsymbol{T}_n=\begin{bmatrix} t(0,0) & & & 0 \\ \vdots & \cdots & \cdots & \\ t(n-1,0) & \cdots t(n-1,k-1) \cdots & & t(n-1,n-1) \end{bmatrix}=\begin{bmatrix} a_0 & & & 0 \\ \vdots & \cdots & \cdots & \\ a_{n-1} & \cdots t(n-1,k-1) \cdots & & b_{n-1} \end{bmatrix} \tag{4}$$

where the $(i,k)$ entries of the matrix and the cells of the triangle are related as:

$$\boldsymbol{T}_n(i,k) = \begin{cases} 0 & (i < k) \\ t(i-1,k-1) & (i \geq k \geq 1) \end{cases} \tag{5}$$

With two boundary sequences and the zero-sum recurrence rule, the explicit formula of the general term to directly compute the interior cells $(i > j > 0)$ is :

$$t(i,j) = \sum_{k=1}^{i-j} \binom{i-k-1}{j-1} (-1)^{i-k} a_k + \sum_{k=1}^{j} \binom{i-k-1}{i-j-1} (-1)^{i-k} b_k \tag{6}$$

which can be easily proved similar to the general term formula for the generalized Pascal's triangle [3], in respect that every row downwards negates the sum of the two cells in the row above.

It is easy to see that the zero-sum triangles are symmetric, $t(i,j) = \pm t(i, i-j)$, if $a_k = \pm b_k$ for all $k$, and the middle cell in an even indexed row is 0 if $a_k + b_k = 0$ and is the largest if $a_k = b_k$ for all $k$.

With the general term formula, a variety of sums of the triangle cells can be found. The sums for the first 2 rows are trivial, and the sum of all the cells in the $i$-th row $S_i$ $(i \geq 2)$ and the sum of all the cells in the whole $n + 1$ row triangle $(n \geq 2)$ $S$ are:

$$\begin{aligned} S_i = \sum_{j=0}^{i} t(i,j) &= a_i + b_i + \sum_{j=1}^{i-1} \left\{ \sum_{k=1}^{i-j} \binom{i-k-1}{j-1} (-1)^{i-k} a_k + \sum_{k=1}^{j} \binom{i-k-1}{i-j-1} (-1)^{i-k} b_k \right\} \\ &= a_i + b_i + (-1)^i \sum_{k=1}^{i-1} \left\{ (-1)^k a_k \sum_{j=1}^{i-k} \binom{i-k-1}{j-1} + (-1)^k b_k \sum_{j=k}^{i-1} \binom{i-k-1}{i-j-1} \right\} \\ &= a_i + b_i + (-1)^i \sum_{k=1}^{i-1} \{(-1)^k a_k 2^{i-k-1} + (-1)^k b_k 2^{i-k-1}\} \\ &= a_i + b_i - \sum_{k=1}^{i-1} (-2)^{i-k-1} \{a_k + b_k\} \end{aligned}$$

$$\begin{aligned} S = \sum_{i=0}^{n} S_i = \sum_{i=0}^{n} \sum_{j=0}^{i} t(i,j) &= a_0 + \sum_{i=1}^{n} \{a_i + b_i\} - \sum_{i=2}^{n} \sum_{k=1}^{i-1} (-2)^{i-k-1} \{a_k + b_k\} \\ &= a_0 + \sum_{i=1}^{n} \{a_i + b_i\} - \sum_{k=1}^{n-1} \{a_k + b_k\} \sum_{i=k+1}^{n} (-2)^{i-k-1} \\ &= a_0 + \sum_{i=1}^{n} \{a_i + b_i\} + \sum_{k=1}^{n-1} \frac{1}{3} ((-2)^{n-k} - 1) \{a_k + b_k\} \\ &= a_0 + a_n + b_n + \frac{1}{3} \sum_{k=1}^{n-1} (2 + (-2)^{n-k}) \{a_k + b_k\} \end{aligned}$$

The two sums are 0 if $a_k + b_k = 0$ for all $k$, and we also have a recursive relation of the sums between two adjacent rows:

$$\begin{aligned} S_i = \sum_{j=0}^{i} t(i,j) = \{a_i + b_i\} - \sum_{k=1}^{i-1} (-2)^{i-k-1} \{a_k + b_k\} &= \{a_i + b_i\} - \{a_{i-1} + b_{i-1}\} + 2 \sum_{k=1}^{i-2} (-2)^{i-k-2} \{a_k + b_k\} \\ &= \{a_i + b_i\} + \{a_{i-1} + b_{i-1}\} - 2 \left\{ \{a_{i-1} + b_{i-1}\} - \sum_{k=1}^{i-2} (-2)^{i-k-2} \{a_k + b_k\} \right\} \\ &= \{a_i + b_i\} + \{a_{i-1} + b_{i-1}\} - 2S_{i-1} \end{aligned}$$

## 3. The Squared Matrices of the Triangles, $\boldsymbol{T}_n^2$

We have our first lemma first:

**Lemma 1**: For a zero-sum triangle $\boldsymbol{T}_n$, the elements of the matrix squared, $\boldsymbol{T}_n^2$, has recursive relations as:

$$\boldsymbol{T}_n^2(i+1, j+1) = \begin{cases} 0 & (i < j) \\ b_i^2 & (i = j) \\ s_0(i) & (i > j = 0) \\ \boldsymbol{T}_n^2(i,j) + t(i,j)\{b_{i-1} + b_i\} - t(i-1,j-1)\{b_{j-1} + b_j\} & (i > j > 0) \end{cases} \tag{7}$$

where we use the following notation for simplification:

$$s_0(i) = \boldsymbol{T}_n^2(i+1,1) = \sum_{k=0}^{i} t(i,k)a_k \tag{8}$$

The details of $s_0(i)$ and more are given in the next section.

**Proof**:

For the first two cases, it is simply true. Since $\boldsymbol{T}_n$ is a lower triangular matrix, $\boldsymbol{T}_n^2$ is also lower triangular, $\boldsymbol{T}_n^2(i+1,j+1) = 0$ if $i < j$, and the diagonal elements are just squared, $\boldsymbol{T}_n^2(i+1,i+1) = t^2(i,i) = b_i^2$.

For all the other elements of $\boldsymbol{T}_n^2$ at $(i+1,j+1)$, $i > j \geq 0$, $i,j = 0,1,2,3,\dots,n-1$,

$$\boldsymbol{T}_n^2(i+1,j+1) = \sum_{k=1}^{n} \boldsymbol{T}_n(i+1,k)\boldsymbol{T}_n(k,j+1) = \sum_{k=j}^{i} t(i,k)t(k,j)$$

For $i > j = 0$, more relations and the proof are given in the next section by using the zero-sum rule.

For $i > j > 0$, we use the zero-sum rule (Eq.2) to substitute the interior elements with their two adjacent ones above or reversely to combine two adjacent ones into one below for simplification:

$$\begin{aligned}
\boldsymbol{T}_n^2(i+1,j+1) &= \sum_{k=j}^{i} t(i,k)t(k,j) = t(i,j)t(j,j) + t(i,i)t(i,j) + \sum_{k=j+1}^{i-1} t(i,k)t(k,j) \\
&= t(i,j)\{b_i + b_j\} + \sum_{k=j+1}^{i-1} \{t(i-1,k-1) + t(i-1,k)\}\{t(k-1,j-1) + t(k-1,j)\} \\
&= t(i,j)\{b_i + b_j\} \\
&+ \sum_{k=j+1}^{i-1} \{t(i-1,k-1)t(k-1,j-1) + t(i-1,k-1)t(k-1,j) - t(i-1,k)t(k,j)\} \\
&= t(i,j)\{b_i + b_j\} + \sum_{k=j+1}^{i-1} t(i-1,k-1)t(k-1,j-1) + t(i-1,j)t(j,j) \\
&- t(i-1,i-1)t(i-1,j) \\
&= t(i,j)b_i + \sum_{k=j}^{i-2} t(i-1,k)t(k,j-1) + \{t(i,j) + t(i-1,j)\}b_j - t(i-1,j)b_{i-1} \\
&= t(i,j)b_i + \sum_{k=j-1}^{i-1} t(i-1,k)t(k,j-1) - t(i-1,j-1)t(j-1,j-1) \\
&- t(i-1,i-1)t(i-1,j-1) + \{t(i,j) + t(i-1,j)\}b_j - t(i-1,j)b_{i-1} \\
&= t(i,j)b_i + \sum_{k=j-1}^{i-1} t(i-1,k)t(k,j-1) - t(i-1,j-1)b_{j-1} \\
&- \{t(i-1,j-1) + t(i-1,j)\}b_{i-1} - t(i-1,j-1)b_j \\
&= t(i,j)b_i + \sum_{k=j-1}^{i-1} t(i-1,k)t(k,j-1) - t(i-1,j-1)b_{j-1} + t(i,j)b_{i-1} \\
&- t(i-1,j-1)b_j = \boldsymbol{T}_n^2(i,j) + t(i,j)\{b_{i-1} + b_i\} - t(i-1,j-1)\{b_{j-1} + b_j\}
\end{aligned}$$

where we use the following identities for the recursive formula:

$$\sum_{k=j-1}^{i-1} t(i-1,k)t(k,j-1) = \boldsymbol{T}_n^2(i,j)$$

■

## 4. The Triangles Multiplied with the Boundary

To make the proof simpler, this section gives another lemma for the triangle matrix multiplied with the left boundary column in three cases: with the boundary itself, $\boldsymbol{s}_0 = \boldsymbol{T}_n\boldsymbol{a}_0$, where $\boldsymbol{a}_0 = [a_0, a_1, \dots, a_{n-1}]^T$, with the boundary shifted up one cell, $\boldsymbol{s}_1 = \boldsymbol{T}_n\boldsymbol{a}_1$, where $\boldsymbol{a}_1 = [a_1, a_2, \dots, a_n]^T$, and with the boundary shifted down one cell, $\boldsymbol{s}_{-1} = \boldsymbol{T}_n\boldsymbol{a}_{-1}$, where $\boldsymbol{a}_{-1} = [0, a_0, a_1, \dots, a_{n-2}]^T$. The elements of the triangle-boundary product vectors can be expressed with the triangle cells as:

$$s_0(i) = \boldsymbol{T}_n^2(i+1,1) = \sum_{k=0}^{i} t(i,k)a_k\,, \qquad s_1(i) = \sum_{k=0}^{i} t(i,k)a_{k+1}\,, \qquad s_{-1}(i) = \sum_{k=1}^{i} t(i,k)a_{k-1} \tag{9}$$

with $s_0(0) = a_0b_0 = a_0^2, s_1(0) = a_1b_0 = a_0a_1, s_{-1}(1) = a_0b_1$, $i \geq 1$.

We can find a few related recursive formulas as in Lemma 2 for the triangle-boundary multiplications, which are more accessible in the proofs in the next sections:

**Lemma 2**:

$$\begin{cases} s_0(0) = a_0^2, \quad s_1(0) = a_0a_1, \quad s_{-1}(1) = a_0b_1 & \text{(base case)} \\ s_{-1}(i) = a_{i-1}\{b_{i-1} + b_i\} - s_{-1}(i-1) - s_0(i-1) & (i \geq 2) \\ s_0(i) = a_0\{a_{i-1} + a_i\} + a_i\{b_{i-1} + b_i\} - s_0(i-1) - s_1(i-1) & (i \geq 1) \\ s_0(i) = a_i\{b_i + b_{i+1}\} - s_{-1}(i) - s_{-1}(i+1) & (i \geq 1) \\ s_1(i) = a_0\{a_i + a_{i+1}\} + a_{i+1}\{b_i + b_{i+1}\} - s_0(i) - s_0(i+1) & (i \geq 1) \end{cases} \tag{10}$$

**Proof**:

We prove them by using the zero-sum rule (Eq.2) to substitute the interior elements with their two adjacent ones if there exist.

$$\begin{aligned} s_{-1}(i) = \sum_{k=1}^{i} t(i,k)a_{k-1} &= t(i,i)a_{i-1} - \sum_{k=1}^{i-1}\{t(i-1,k-1) + t(i-1,k)\}a_{k-1} \\ &= a_{i-1}b_i - \sum_{k=0}^{i-1} t(i-1,k)a_k + t(i-1,i-1)a_{i-1} - \sum_{k=1}^{i-1} t(i-1,k)a_{k-1} \\ &= a_{i-1}\{b_{i-1} + b_i\} - s_0(i-1) - s_{-1}(i-1) \end{aligned}$$

$$\begin{aligned} s_0(i) = \sum_{k=0}^{i} t(i,k)a_k &= t(i,0)a_0 + t(i,i)a_i - \sum_{k=1}^{i-1}\{t(i-1,k-1) + t(i-1,k)\}a_k \\ &= a_ia_0 + a_ib_i - \sum_{k=0}^{i-1} t(i-1,k)a_{k+1} + t(i-1,i-1)a_i - \sum_{k=0}^{i-1} t(i-1,k)a_k + t(i-1,0)a_0 \\ &= a_0\{a_{i-1} + a_i\} + a_i\{b_{i-1} + b_i\} - s_1(i-1) - s_0(i-1) \end{aligned}$$

$$\begin{aligned} s_0(i) = \sum_{k=0}^{i} t(i,k)a_k &= t(i,i)a_i - \sum_{k=0}^{i-1}\{t(i,k+1) + t(i+1,k+1)\}a_k \\ &= a_ib_i - \sum_{k=1}^{i} t(i,k)a_{k-1} - \sum_{k=1}^{i+1} t(i+1,k)a_{k-1} + t(i+1,i+1)a_i \\ &= a_i\{b_i + b_{i+1}\} - s_{-1}(i) - s_{-1}(i+1) \end{aligned}$$

$$\begin{aligned} s_1(i) = \sum_{k=0}^{i} t(i,k)a_{k+1} &= t(i,i)a_{i+1} - \sum_{k=0}^{i-1}\{t(i,k+1) + t(i+1,k+1)\}a_{k+1} \\ &= a_{i+1}b_i - \sum_{k=0}^{i} t(i,k)a_k + t(i,0)a_0 - \sum_{k=0}^{i+1} t(i+1,k)a_k + t(i+1,0)a_0 + t(i+1,i+1)a_{i+1} \\ &= a_0\{a_i + a_{i+1}\} + a_{i+1}\{b_i + b_{i+1}\} - s_0(i) - s_0(i+1) \end{aligned}$$

■

## 5. The Triangles for Idempotent Matrices

A matrix $\boldsymbol{T}$ is idempotent if $\boldsymbol{T}^2 = \boldsymbol{T}$. In the squared triangle matrix as in Eq.(7), for the diagonal entries, to make $b_i^2 = b_i$, we have $b_i = 0$ or $b_i = 1$.

For $i > j > 0$, if we have $b_{i-1} + b_i = b_{j-1} + b_j = 1$, then we have $\boldsymbol{T}_n^2(i+1,j+1) = \boldsymbol{T}_n^2(i,j) + t(i,j) - t(i-1,j-1)$. To make it, we can simply set the integer sequences $b_i$ as alternatively 0 and 1 :

$$t(i,i) = b_i = 1 - b_{i-1} = \mathrm{mod}(i + a_0, 2) \tag{11}$$

with which we also have $b_ib_{i-1} = 0$.

With the right boundary edge numbers given in Eq.(11), the elements of the squared triangle matrix, Eq.(7), are:

$$\boldsymbol{T}_n^2(i+1,j+1)=\begin{cases}0 & (i<j)\\ b_i & (i=j)\\ s_0(i) & (i>j=0)\\ \boldsymbol{T}_n^2(i,j)+t(i,j)-t(i-1,j-1) & (i>j>0)\end{cases}$$

For $i>j\geq 0$, it is easy to verify with Eq.(7&10) if $s_0(i)=a_i$ and $\boldsymbol{T}_n^2(i,j)=t(i-1,j-1)$, and we then have $\boldsymbol{T}_n^2=\boldsymbol{T}_n$, which is an idempotent matrix. To meet the conditions, we use $s_1(i)$ for the solution and then we prove $s_{-1}(i)=0, s_0(i)=a_i, s_1(i)=a_0\{a_i+a_{i+1}\}-a_i$, and $\boldsymbol{T}_n^2=\boldsymbol{T}_n$ as follows.

If $a_0=0$, $t(2m,2m)=b_{2m}=\text{mod}(2m+a_0,2)=a_0=0$, then we use equation $s_0(2m)=a_{2m}$ for the solution of the even indexed boundary numbers $a_{2m}$ :

$$a_{2m}=s_0(2m)=\sum_{k=1}^{2m-1}t(2m,k)a_k \tag{12}$$

In such a case, $b_{2m+1}=\text{mod}(2m+1+a_0,2)=1$, we need to prove the following equalities hold:

$$s_{-1}(2m+1)=0, s_0(2m+1)=a_{2m+1}, s_1(2m)=-a_{2m}, \boldsymbol{T}_n^2(2m+2,j+1)=t(2m+1,j) \tag{13}$$

If $a_0=1$, $t(2m-1,2m-1)=b_{2m-1}=\text{mod}(2m-1+a_0,2)=1-a_0=0$, then we use equation $s_1(2m-1)=a_{2m}$ for the solution of $a_{2m}$ :

$$a_{2m}=s_1(2m-1)=\sum_{k=0}^{2m-2}t(2m-1,k)a_{k+1} \tag{14}$$

In such a case, $b_{2m}=\text{mod}(2m+a_0,2)=1$, we need the following equalities hold:

$$s_{-1}(2m)=0, s_0(2m)=a_{2m}, s_1(2m)=a_{2m+1}, \boldsymbol{T}_n^2(2m+1,j+1)=t(2m,j) \tag{15}$$

By induction we give the proof below that with Eq.(12&14) no matter what $a_{2m+1}$ are, or $a_{2m+1}$ can be any independent input, the equalities in Eq.(13) or Eq.(15) hold.

In summary, a sufficient condition for the triangle matrix to be idempotent is the boundary integer sequences to be given by Eq.(11&12), or by Eq.(11&14), as descripted in the following theorem with Eq.(12&14) combined into one Eq.(16) and the boundary relations Eq.(13&15) combined into Eq.(17).

**Theorem 1**: The matrices of the zero-sum triangles are idempotent if the triangles are generated with the boundary sequences:

$$a_{2m}=\begin{cases}s_0(2m)=\displaystyle\sum_{k=1}^{2m-1}t(2m,k)a_k & (b_{2m}=a_0=0)\\ s_1(2m-1)=\displaystyle\sum_{k=0}^{2m-2}t(2m-1,k)a_{k+1} & (b_{2m}=a_0=1)\end{cases} \tag{16}$$

where $m\geq 1$, $a_{2m+1}$ can be freely chosen, and $b_{2m+1}=1-b_{2m}$.

The multiplication of the triangle matrix with the left boundary gives:

$$\begin{cases}s_0(0)=a_0, & s_1(0)=a_0a_1, & s_{-1}(1)=0 & \text{(base case)}\\ s_0(i)=a_i, & s_1(i)=a_0\{a_i+a_{i+1}\}-a_i, & s_{-1}(i)=0 & (i\geq 1)\end{cases} \tag{17}$$

**Proof**:

For $i>j>0$, staring with $\boldsymbol{T}_n^2(1,1)=a_0^2=a_0=t(0,0)=\boldsymbol{T}_n(1,1)$ and the recursive formula, $\boldsymbol{T}_n^2(i+1,j+1)=\boldsymbol{T}_n^2(i,j)+t(i,j)-t(i-1,j-1)$, by induction we can easily have $\boldsymbol{T}_n^2(i+1,j+1)=t(i,j)$ from $\boldsymbol{T}_n^2(i,j)=t(i-1,j-1)$. So we only need to prove for $s_0(i)=a_i$ such that the left boundary elements of the squared triangle also meet the idempotent property: $\boldsymbol{T}_n^2(i+1,1)=s_0(i)=a_i=t(i,0)=\boldsymbol{T}_n(i+1,1)$.

For $n=1,2,3$, the triangles are formed as:

$$\boldsymbol{T}_1=a_0,\quad \boldsymbol{T}_2=\begin{bmatrix}a_0 & 0\\ a_1 & b_1\end{bmatrix},\quad \boldsymbol{T}_3=\begin{bmatrix}a_0 & 0 & 0\\ a_1 & b_1 & 0\\ a_2 & -a_1-b_1 & b_2\end{bmatrix}$$

Let $a_1$ be free, it is easy to verify that $\boldsymbol{T}_1^2=\boldsymbol{T}_1$, $\boldsymbol{T}_2^2=\boldsymbol{T}_2$, $s_0(0)=a_0^2=a_0, s_0(1)=a_1a_0+b_1a_1=a_1, s_{-1}(1)=a_0b_1=0, s_1(0)=a_0a_1$, and with the recursive formulas, $s_{-1}(2)=a_1-s_{-1}(1)-s_0(1)=0$.

With $m=1, n=2m+1=3$, we can have $\boldsymbol{T}_3^2=\boldsymbol{T}_3$ if we set as in Eq.(16),

$$a_2=\begin{cases}s_0(2)=\displaystyle\sum_{k=1}^{1}t(2,k)a_k=t(2,1)a_1=-a_1(a_1+b_1) & (b_2=a_0=0)\\ s_1(1)=\displaystyle\sum_{k=0}^{0}t(1,k)a_{k+1}=t(1,0)a_1=a_1^2 & (b_2=a_0=1)\end{cases}$$

Using definitions Eq.(9) and the recursive formulas Eq.(10), we have for $a_0 = 0$, $s_1(1) = a_0\{a_1 + a_2\} + a_2 - s_0(1) - s_0(2) = a_0\{a_1 + a_2\} - a_1$ , $s_{-1}(3) = a_2 - s_{-1}(2) - s_0(2) = a_2 - 0 - a_2 = 0$ , and for $a_0 = 1$ , $s_0(2) = a_0\{a_1 + a_2\} + a_2 - s_0(1) - s_1(1) = a_2 = a_0\{a_1 + a_2\} - a_1$, $s_{-1}(3) = a_2 - s_{-1}(2) - s_0(2) = 0$. Or by combing the above we have $s_1(1) = a_0\{a_1 + a_2\} - a_1$, $s_0(1) = a_1$, $s_{-1}(1) = 0$. Eq.(17) holds for $i = 1$.

For $m = 1$, $n = 2m + 2 = 4$, with the zero-sum rule we give $\boldsymbol{T}_4$ as:

$$\boldsymbol{T}_4 = \begin{bmatrix} a_0 & 0 & 0 & 0 \\ a_1 & b_1 & 0 & 0 \\ a_2 & -a_1 - b_1 & b_2 & 0 \\ a_3 & a_1 + b_1 - a_2 & a_1 + b_1 - b_2 & b_3 \end{bmatrix}$$

We have $s_1(2) = a_2a_1 + t(2,1)a_2 + b_2a_3 = a_2(-b_1) + a_0a_3 = a_0\{a_2 + a_3\} - a_2$, $s_0(3) = a_0\{a_2 + a_3\} + a_3 - s_0(2) - s_1(2) = a_3$ and $s_{-1}(4) = a_3 - s_{-1}(3) - s_0(3) = 0$.

Put equalities of index 2 together, and we see Eq.(17) holds for $i = 2$ as $s_1(2) = a_0\{a_2 + a_3\} - a_2$, $s_0(2) = a_2$, $s_{-1}(2) = 0$.

With Eq.(16), we find the solution for the next even indexed boundary number:

$$a_{2m} = \begin{cases} s_0(2m) & (a_0 = 0) \\ s_1(2m-1) & (a_0 = 1) \end{cases}$$

and we assume two steps of the three equalities in Eq.(17) hold for integer $i = 2m - 1$ and $i = 2m$ $(m > 1)$:

$$s_1(2m-1) = a_0\{a_{2m-1} + a_{2m}\} - a_{2m-1}, \qquad s_0(2m-1) = a_{2m-1}, \qquad s_{-1}(2m-1) = 0$$
$$s_1(2m) = a_0\{a_{2m} + a_{2m+1}\} - a_{2m}, \qquad s_0(2m) = a_{2m}, \qquad s_{-1}(2m) = 0$$

Then with the recursive relations in Eq.(10), we can have $s_{-1}(i)$ and $s_0(i)$ in Eq.(17) hold for $i = 2m + 1$ as :

$$s_{-1}(2m+1) = a_{2m}\{b_{2m} + b_{2m+1}\} - s_{-1}(2m) - s_0(2m) = a_{2m} - 0 - a_{2m} = 0$$
$$\begin{aligned} s_0(2m+1) &= a_0\{a_{2m} + a_{2m+1}\} + a_{2m+1}\{b_{2m} + b_{2m+1}\} - s_0(2m) - s_1(2m) \\ &= a_0\{a_{2m} + a_{2m+1}\} + a_{2m+1} - a_{2m} - a_0\{a_{2m} + a_{2m+1}\} + a_{2m} = a_{2m+1} \end{aligned}$$

Hence by induction, we have $s_{-1}(i) = 0$ and $s_0(i) = a_i$ hold. Then with the recursive relation for $s_1(i)$ we have the last one in Eq.(17) :

$$\begin{aligned} s_1(2m+1) &= a_0\{a_{2m+1} + a_{2m+2}\} + a_{2m+2}\{b_{2m} + b_{2m+1}\} - s_0(2m+1) - s_0(2m+2) \\ &= a_0\{a_{2m+1} + a_{2m+2}\} + a_{2m+2} - a_{2m+1} - a_{2m+2} = a_0\{a_{2m+1} + a_{2m+2}\} - a_{2m+1} \end{aligned}$$

which shows that $s_1(i) = a_0\{a_i + a_{i+1}\} - a_i$ also holds.

■

With Theorem 1, for $a_0 = 1$, we can find the first few dependent even indexed boundary cells are: $a_2 = a_1^2, a_4 = (1 - a_1)a_1^3 + (2a_1 - 1)a_3, a_6 = 4a_1^4 - 5a_1^5 + 2a_1^6 + a_3 + 7a_1^2a_3 + a_3^2 - a_1^3(4a_3 + 1) - 2a_5 + 2a_1(a_5 - 2a_3), a_8 = -39a_1^6 + 22a_1^7 - 5a_1^8 + a_1^5(12a_3 + 35) - a_1^4(41a_3 + 16) + a_1^3(56a_3 - 4a_5 + 3) + 5a_5 + a_3(2a_5 - 5a_3 - 3) + a_1^2(11a_5 - 2a_3(20 + 3a_3)) - 3a_7 + a_1(a_3(11a_3 + 16) + 2(a_7 - 6a_5))$.

For $a_0 = 0$ , the first few even indexed boundary cells are: $a_2 = -a_1(1 + a_1), a_4 = a_1(1 + a_1)^3 - 2(1 + a_1)a_3, a_6 = -24a_1^4 - 11a_1^5 - 2a_1^6 - (-5 + a_3)a_3 + a_1^3(-26 + 4a_3) + a_1^2(-14 + 13a_3) + a_1(-3 + 14a_3 - 2a_5) - 3a_5, a_8 = 151a_1^6 + 42a_1^7 + 5a_1^8 + a_1^4(359 - 71a_3) + a_1^5(301 - 12a_3) + 12a_3^2 + 14a_5 - 2a_3(14 + a_5) + 4a_1^3(64 - 42a_3 + a_5) + a_1^2(101 + a_3(-199 + 6a_3) + 17a_5) + a_1(17 + a_3(-118 + 17a_3) + 26a_5 - 2a_7) - 4a_7$.

Since no division is involved in Eq.(16), the triangle only has integers as long as the odd indexed boundary numbers are input with integers. An example of idempotent triangle is Triangle 1 shown on the left in Section 2, where we use $a_i = 1$. More examples can be with $a_i = 0$, $a_i = (-1)^i$, or with $a_i$ using a repeating pattern of 4 integers $\{0, -1, 0, 1\}$. The first 9 rows of two more triangles (3&4) are given below:

```
                1
             -1   0
            1   1   1
         -1  -2  -2   0
        1   3   4   2   1
     -1  -4  -7  -6  -3   0
    1   5  11  13   9   3   1
 -1  -6 -16 -24 -22 -12  -4   0
1   7  22  40  46  34  16   4   1
```

Triangle 3

```
                0
             -1   1
            0   0   0
          1   0   0   1
        0  -1   0  -1   0
     -1   1   1   1   1   1
    0   0  -2  -2  -2  -2   0
  1   0   2   4   4   4   2   1
0  -1  -2  -6  -8  -8  -6  -3   0
```

Triangle 4

In fact, if we remove the leftmost boundary cells, the triangle with $a_0 = 1$ becomes the triangle with $a_0 = 0$, and vice versa.

Based on Theorem 1, we can also see that the even indexed rows are linear combinations of the odd indexed rows above, and the even indexed columns are linear combinations of the odd indexed columns rightwards.

**6. The Triangles for Involutory Matrices**

For a given idempotent matrix $\boldsymbol{T}$, $\boldsymbol{T}^2 = \boldsymbol{T}$, we can easily find two involutory matrices, $2\boldsymbol{T} - \boldsymbol{I}$ and $\boldsymbol{I} - 2\boldsymbol{T}$, but we can directly find more general involutory matrices with the triangles constructed with the zero-sum rule.

Based on the recursive formula of the squared triangle matrix, Eq.(7), for the diagonal entries, to make $b_i^2 = 1$, we have $b_i = 1$ or $b_i = -1$. For $i > j > 0$, if we have $b_{i-1} + b_i = b_{j-1} + b_j = 0$, then we have $\boldsymbol{T}_n^2(i+1, j+1) = \boldsymbol{T}_n^2(i,j)$. To make $b_{i-1} + b_i = b_{j-1} + b_j = 0$, we can simply set the boundary sequence $b_i$ as alternatively $-1$ and 1,

$$t(i,i) = b_i = (-1)^i b_0 = (-1)^i a_0 \tag{18}$$

Then the elements of the squared triangle matrix are:

$$\boldsymbol{T}_n^2(i+1, j+1) = \begin{cases} 0 & (i < j) \\ 1 & (i = j) \\ s_0(i) & (i > j = 0) \\ \boldsymbol{T}_n^2(i,j) & (i > j > 0) \end{cases}$$

For $i > j$, if $s_0(i) = 0$ and $\boldsymbol{T}_n^2(i,j) = 0$, we have $\boldsymbol{T}_n^2 = \boldsymbol{I}$, which is an involutory matrix.

To meet the conditions, we use $s_0(2m) = 0$ for the solution and then we prove $s_0(2m+1) = 0$, $s_{-1}(i) = (-1)^i$, $s_1(i) = a_0\{a_i + a_{i+1}\}$, and $\boldsymbol{T}_n^2 = \boldsymbol{I}$ as follows.

If $i = 2m$ is even, $b_i = b_{2m} = t(2m, 2m) = (-1)^{2m} b_0 = a_0$, we can solve the equation $s_0(2m) = 0$ to find the even indexed boundary numbers $a_{2m}$:

$$s_0(2m) = 0 = \sum_{k=0}^{2m} t(2m,k)a_k = t(2m,0)a_0 + t(2m,2m)a_{2m} + \sum_{k=1}^{2m-1} t(2m,k)a_k$$

$$= 2a_0 a_{2m} + \sum_{k=1}^{2m-1} t(2m,k)a_k$$

With which we easily give the solution :

$$a_{2m} = -\frac{a_0}{2} \sum_{k=1}^{2m-1} t(2m,k)a_k \tag{19}$$

In summary, a sufficient condition for the triangle matrix to be involutory is the boundary integer sequences to be given by Eq.(18&19), as descripted in the following theorem.

**Theorem 2**: The matrices of the triangles are involutory if the triangles are generated with the boundary sequences

$$a_{2m} = -\frac{a_0}{2} \sum_{k=1}^{2m-1} t(2m,k)a_k, \qquad b_{2m} = a_0 \tag{20}$$

where $m \geq 1$, $a_0 = \pm 1$, $a_{2m+1}$ can be freely chosen, and $b_{2m+1} = 1 - b_{2m}$.

The multiplication of the triangle matrix with the left boundary gives:

$$\begin{cases} s_0(0) = 1, \quad s_1(0) = a_0 a_1, \quad s_{-1}(1) = -1 & \text{(base case)} \\ s_0(i) = 0, \quad s_1(i) = a_0\{a_i + a_{i+1}\}, \quad s_{-1}(i) = (-1)^i & (i \geq 1) \end{cases} \tag{21}$$

**Proof**:

For $i > j > 0$, if we have $s_0(i-1) = \boldsymbol{T}_n^2(i, 1) = 0$, with the recursive formula, $\boldsymbol{T}_n^2(i+1, j+1) = \boldsymbol{T}_n^2(i,j)$, by induction we can easily have $\boldsymbol{T}_n^2(i,j) = 0$ and $\boldsymbol{T}_n^2(i+1, j+1) = 0$ for all $i > j > 1$. So we only need to prove with the left boundary elements $s_0(i) = 0$ for $i > 0$.

For even indexed boundary elements, $a_{2m}$, we use Eq.(20) for the solution, so its corresponding equation $s_0(2m) = 0$ always holds. Then the proof is needed for the odd indexed boundary elements, or $s_0(2m+1) = 0$, and we do it by induction.

For $n = 1,2,3$, with $a_0 = \pm 1$, the triangles are formed as:

$$\boldsymbol{T}_1 = a_0, \qquad \boldsymbol{T}_2 = \begin{bmatrix} a_0 & 0 \\ a_1 & -a_0 \end{bmatrix}, \qquad \boldsymbol{T}_3 = \begin{bmatrix} a_0 & 0 & 0 \\ a_1 & -a_0 & 0 \\ a_2 & a_0 - a_1 & a_0 \end{bmatrix}$$

Let $a_1$ be free, it is easy to verify that $\boldsymbol{T}_1^2 = \boldsymbol{I}$ and $\boldsymbol{T}_2^2 = \boldsymbol{I}$, $s_0(0) = a_0^2 = 1, s_1(0) = a_0 a_1, s_{-1}(1) = a_0 b_1 =$

$-1, s_0(1) = a_1a_0 + b_1a_1 = 0$.

With $m = 1, n = 2m + 1 = 3$, we can have $\boldsymbol{T}_3^2 = \boldsymbol{I}$ if we find $a_2$ by using Eq.(20):

$$a_2 = a_{2m} = -\frac{a_0}{2}\sum_{k=1}^{2m-1} t(2m,k)a_k = -\frac{a_0}{2}t(2,1)a_1 = \frac{a_0}{2}a_1(a_1 - a_0)$$

Then we can verify $s_1(1) = a_1^2 + b_1a_2 = a_0\{a_1 + a_2\}$, $s_0(2) = a_2a_0 + t(2,1)a_1 + b_2a_2 = 0$, $s_{-1}(2) = t(2,1)a_0 + b_2a_1 = 1$, and $\boldsymbol{T}_3^2 = \boldsymbol{I}$, and give $\boldsymbol{T}_4$ as:

$$\boldsymbol{T}_4 = \begin{bmatrix} a_0 & 0 & 0 & 0 \\ a_1 & -a_0 & 0 & 0 \\ \frac{a_0}{2}a_1(a_1 - a_0) & a_0 - a_1 & a_0 & 0 \\ a_3 & -\frac{a_0}{2}(a_1 - a_0)(a_1 - 2a_0) & a_1 - 2a_0 & -a_0 \end{bmatrix}$$

It is also easy to verify $s_1(2) = a_0\{a_2 + a_3\}$, $s_0(3) = 0$, $s_{-1}(3) = -1$, and $\boldsymbol{T}_4^2 = \boldsymbol{I}$.

For the boundary conditions, we can summarize the equalities above as :

$$s_0(1) = 0, \quad s_1(1) = a_0\{a_1 + a_2\}, \quad s_{-1}(1) = -1$$
$$s_0(2) = 0, \quad s_1(2) = a_0\{a_2 + a_3\}, \quad s_{-1}(2) = 1$$

If we use Eq.(20) to find $a_{2m}$, we have $s_0(2m) = 0$. Assume Eq.(21) holds for integer $i = 2m$,

$$s_0(2m) = 0, \quad s_1(2m) = a_0\{a_{2m} + a_{2m+1}\}, \quad s_{-1}(2m) = (-1)^{2m} = 1$$

Then we can find those for $i = 2m + 1$ with the recursive formulas in Eq.(10). For $s_{-1}(i)$ we have:

$$s_{-1}(2m + 1) = a_{2m}\{b_{2m} + b_{2m+1}\} - s_{-1}(2m) - s_0(2m) = -1$$

By induction, we have $s_{-1}(i) = (-1)^i$ holds, and with Eq.(10) again for $s_0(i)$ we have :

$$s_0(i) = a_i\{b_i + b_{i+1}\} - s_{-1}(i) - s_{-1}(i + 1) = 0$$

Finally, with the recursive formula Eq.(10) for $s_1(i)$, we have

$$s_1(i) = a_0\{a_i + a_{i+1}\} + a_{i+1}\{b_i + b_{i+1}\} - s_0(i) - s_0(i + 1) = a_0\{a_i + a_{i+1}\}$$

■

With $a_0 = \pm 1$, the first few dependent boundary cells are: $a_2 = \frac{1}{2}a_1(a_0a_1 - 1), a_4 = \frac{1}{8}\big(a_1\big(2 + 4a_1^2 - a_0(5a_1 + a_1^3 - 8a_3)\big) - 12a_3\big), a_6 = \frac{1}{16}(-8a_1 - 38a_1^3 - 8a_1^5 + 40a_3 + 40a_1^2a_3 - 40a_5 + a_0(28a_1^2 + 25a_1^4 + a_1^6 - 8a_1^3a_3 + 8a_3^2 + a_1(-68a_3 + 16a_5)))$.

Since a division is involved in Eq.(20), if all the odd indexed boundary numbers are integers, the triangle does not have to be all integers. If we expand the solution for $a_{2m}$, Eq.(20), to find all the terms involving $a_{2m-1}$, then we can choose some $a_{2m-1}$ to make $a_{2m}$ integers:

$$a_{2m} = -\frac{a_0}{2}\sum_{k=1}^{2m-1} t(2m,k)a_k = \frac{a_0}{2}\left\{-t(2m,1)a_1 - t(2m,2m-1)a_{2m-1} - \sum_{k=2}^{2m-2} t(2m,k)a_k\right\}$$
$$= \frac{a_0}{2}\left\{\{t(2m-1,0) + t(2m-1,1)\}a_1 - t(2m,2m-1)a_{2m-1} - \sum_{k=2}^{2m-2} t(2m,k)a_k\right\}$$
$$= \frac{a_0}{2}\left\{\{a_1 - t(2m,2m-1)\}a_{2m-1} + a_1t(2m-1,1) - \sum_{k=2}^{2m-2} t(2m,k)a_k\right\}$$

With the general formula for the interior cells Eq.(6) we can find $t(2m, 2m - 1)$ as :

$$t(2m,2m-1) = \sum_{k=1}^{1}\binom{2m-k-1}{2m-1-1}(-1)^{2m-k}a_k + \sum_{k=1}^{2m-1}\binom{2m-k-1}{2m-2m+1-1}(-1)^{2m-k}b_k$$
$$= -a_1 + \sum_{k=1}^{2m-1}(-1)^{2k}a_0 = a_0(2m-1) - a_1$$

Then all the terms involving $a_{2m-1}$ can be

$$\{a_1 - t(2m,2m-1)\}a_{2m-1} = \{2a_1 - a_0(2m-1)\}a_{2m-1} = 2(a_1 - a_0m)a_{2m-1} + a_0a_{2m-1}$$

which gives

$$a_{2m} = \frac{a_0}{2}\left\{a_0a_{2m-1} + 2(a_1 - a_0m)a_{2m-1} + a_1t(2m-1,1) - \sum_{k=2}^{2m-2} t(2m,k)a_k\right\}$$

We can just choose an even or an odd number for $a_{2m-1}$ to make the number inside the curly brackets an even number. Therefore, to make all the triangle cells integers, the odd indexed left boundary cells are not completely free to choose, but just limited to odd or even numbers.

An example of involutory triangle is Triangle 2 shown on the right in Section 2, where we use $a_i = 1 + \mathrm{mod}(i, 2)$. More examples can be with $a_i = 0$, $a_i = 3$, $a_i = (-1)^i$, or with $a_i$ using a repeating pattern of 4 integers $\{1,0,-2,1\}$. The first 9 rows of two more triangles (5&6) are given below:

```
                1
              3   -1
            3   -2   1
          3   -1   1   -1
        3   -2   0   0   1
      3   -1   2   0   -1  -1
    3   -2  -1  -2   1   2   1
  3   -1   3   3   1  -3  -3  -1
3   -2  -2  -6  -4   2   6   4   1
```

Triangle 5

```
                1
              1   -1
            0   0   1
         -2   0   -1  -1
        1   2   1   2   1
      1   -3  -3  -3  -3  -1
    0   2   6   6   6   4   1
 -2  -2  -8 -12 -12 -10  -5  -1
1   4  10  20  24  22  15   6   1
```

Triangle 6

Similar to the idempotent triangles, if we remove the leftmost boundary cells, the triangle with $a_0 = 1$ becomes the triangle with $a_0 = -1$, and vice versa. In addition, if we negate all the numbers in an involutory triangle, it is also involutory, or if $\boldsymbol{T}^2 = \boldsymbol{I}$, $(-\boldsymbol{T})^2 = \boldsymbol{I}$. This actually switches between the triangles with $a_0 = 1$ and the triangles with $a_0 = -1$.

## 7. The Triangles for Nilpotent Matrices of Index 2

If a triangular matrix is nilpotent, it should be strictly triangular and all the main diagonal entries should be 0. The triangles we have in the above sections are not, but we can shift a triangular matrix downwards or leftwards such that they become strictly triangular with all the main diagonal entries being 0.

For square matrices, to shift a matrix down by one row, we fill a row of zeros at the top and remove the bottom row. To shift left by one column, we fill a column of zeros on the right and remove the first column. The shift operator is the same as matrix multiplication with the shift-left matrix, $\boldsymbol{B}_L$, a lower subdiagonal matrix with 1s as all subdiagonal entries, left-multiply to shift-down and right-multiply to shift-left. With the idempotent triangle matrices, we can convert them into nilpotent matrices by shifting: $\boldsymbol{N}_D = \boldsymbol{B}_L\boldsymbol{T}_n$ and $\boldsymbol{N}_L = \boldsymbol{T}_n\boldsymbol{B}_L$.

With the idempotent triangles given by Theorem 1 above, we have the following theorem.

**Theorem 3**: If an idempotent matrix of zero-sum triangle $\boldsymbol{T}_n$ is shifted one cell downwards, $\boldsymbol{N}_D = \boldsymbol{B}_L\boldsymbol{T}_n$, or one cell leftwards, $\boldsymbol{N}_L = \boldsymbol{T}_n\boldsymbol{B}_L$, the matrices are nilpotent matrices of index 2, $\boldsymbol{N}_D^2 = \boldsymbol{N}_L^2 = 0$.

**Proof**:

Let $\boldsymbol{T}_n$ be an idempotent matrix of the triangle generated with the zero-sum rule, and we still use the same notations as in Section 5, cells $\boldsymbol{T}_n(i+1, j+1) = t(i, j)$ and the boundary cells meet the conditions as in Theorem 1. Use $\boldsymbol{N}_D$ and $\boldsymbol{N}_L$ to represent the matrices by shifting $\boldsymbol{T}_n$ down one row and shifting $\boldsymbol{T}_n$ left one column, respectively, which are equivalent to $\boldsymbol{N}_D = \boldsymbol{B}_L\boldsymbol{T}_n$ and $\boldsymbol{N}_L = \boldsymbol{T}_n\boldsymbol{B}_L$. The elements of the matrices are:

$$\boldsymbol{N}_D(i+1, j+1) = \begin{cases} 0 & (i = 0) \\ \boldsymbol{T}_n(i, j+1) & (i > 0) \end{cases}$$

$$\boldsymbol{N}_L(i+1, j+1) = \begin{cases} 0 & (j = n-1) \\ \boldsymbol{T}_n(i+1, j+2) & (j < n-1) \end{cases}$$

The elements of $\boldsymbol{T}_n\boldsymbol{B}_L\boldsymbol{T}_n$ for $i > j$ are

$$(\boldsymbol{T}_n\boldsymbol{B}_L\boldsymbol{T}_n)(i+1, j+1) = (\boldsymbol{T}_n\boldsymbol{N}_D)(i+1, j+1) = \sum_{k=0}^{n-1} \boldsymbol{T}_n(i+1, k+1)\boldsymbol{N}_D(k+1, j+1)$$

$$= \sum_{k=1}^{n-1} \boldsymbol{T}_n(i+1, k+1)\boldsymbol{T}_n(k, j+1) = \sum_{k=j+1}^{i} t(i, k)t(k-1, j)$$

For $j = 0$, $t(k-1, 0) = a_{i-1}$, with the boundary multiplication equalities in Theorem 1, we have the first column all 0s :

$$(\boldsymbol{T}_n\boldsymbol{B}_L\boldsymbol{T}_n)(i+1, 1) = \sum_{k=1}^{i} t(i, k)t(k-1, 0) = \sum_{k=1}^{i} t(i, k)a_{i-1} = s_{-1}(i) = 0$$

For $j > 0$, we have all the elements in a column equal to the column on the left up to a sign:

$$(\boldsymbol{T}_n\boldsymbol{B}_L\boldsymbol{T}_n)(i+1,j+1)=\sum_{k=j+1}^{i} t(i,k)t(k-1,j)=-\sum_{k=j+1}^{i} t(i,k)\{t(k-1,j-1)+t(k,j)\}$$
$$=t(i,j)t(j-1,j-1)-\sum_{k=j}^{i} t(i,k)t(k-1,j-1)+t(i,j)t(j,j)-\sum_{k=j}^{i} t(i,k)t(k,j)$$
$$=t(i,j)\{b_{j-1}+b_j\}-(\boldsymbol{T}_n\boldsymbol{B}_L\boldsymbol{T}_n)(i+1,j)-\boldsymbol{T}_n^2(i+1,j+1)$$
$$=t(i,j)\{b_{j-1}+b_j\}-(\boldsymbol{T}_n\boldsymbol{B}_L\boldsymbol{T}_n)(i+1,j)-\boldsymbol{T}_n(i+1,j+1)$$
$$=t(i,j)\{b_{j-1}+b_j\}-(\boldsymbol{T}_n\boldsymbol{B}_L\boldsymbol{T}_n)(i+1,j)-t(i,j)=-(\boldsymbol{T}_n\boldsymbol{B}_L\boldsymbol{T}_n)(i+1,j)$$

With the first column all 0s, this gives $\boldsymbol{T}_n\boldsymbol{B}_L\boldsymbol{T}_n=0$.

Thus, we have $\boldsymbol{N}_D^2=\boldsymbol{B}_L\boldsymbol{T}_n\boldsymbol{B}_L\boldsymbol{T}_n=0$ and $\boldsymbol{N}_L^2=\boldsymbol{T}_n\boldsymbol{B}_L\boldsymbol{T}_n\boldsymbol{B}_L=0$, and both $\boldsymbol{N}_D$ and $\boldsymbol{N}_L$ are nilpotent matrices of index 2.

■

## 8. The Triangles for Unipotent Matrices of Index 2

With the nilpotent matrices generated by shifting an idempotent triangle, we can have their corresponding unipotent matrices of index 2, $\boldsymbol{S}_D=\boldsymbol{I}+\boldsymbol{N}_D=\boldsymbol{I}+\boldsymbol{B}_L\boldsymbol{T}_n$ and $\boldsymbol{S}_L=\boldsymbol{I}+\boldsymbol{N}_L=\boldsymbol{I}+\boldsymbol{T}_n\boldsymbol{B}_L$. Then we have the following theorem for a bidiagonal matrix with both diagonals all 1s, a full-sized Jordan block of eigenvalue 1.

**Theorem 4**: If $\boldsymbol{N}_D$ and $\boldsymbol{N}_L$ are nilpotent matrices generated by shifting an idempotent matrix of zero-sum triangle one cell downwards and one cell leftwards, respectively, the product of the corresponding unipotent matrices makes a bidiagonal matrix with both diagonals all 1s: $\boldsymbol{S}_D\boldsymbol{S}_L=(\boldsymbol{I}+\boldsymbol{B}_L\boldsymbol{T}_n)(\boldsymbol{I}+\boldsymbol{T}_n\boldsymbol{B}_L)=\boldsymbol{I}+\boldsymbol{B}_L$.

**Proof**:

With the idempotent property, we can easily have:

$$\boldsymbol{S}_D\boldsymbol{S}_L=(\boldsymbol{I}+\boldsymbol{B}_L\boldsymbol{T}_n)(\boldsymbol{I}+\boldsymbol{T}_n\boldsymbol{B}_L)=\boldsymbol{I}+\boldsymbol{B}_L\boldsymbol{T}_n+\boldsymbol{T}_n\boldsymbol{B}_L+\boldsymbol{B}_L\boldsymbol{T}_n\boldsymbol{T}_n\boldsymbol{B}_L=\boldsymbol{I}+\boldsymbol{B}_L\boldsymbol{T}_n+\boldsymbol{T}_n\boldsymbol{B}_L+\boldsymbol{B}_L\boldsymbol{T}_n\boldsymbol{B}_L$$

They are all lower triangular matrices, and the main diagonal elements of the triangular unipotent matrices are all 1s, so we only consider the lower elements for the proof.

For the elements below the lower subdiagonal, $i-2\geq j\geq 0$, we have :

$$(\boldsymbol{B}_L\boldsymbol{T}_n+\boldsymbol{T}_n\boldsymbol{B}_L+\boldsymbol{B}_L\boldsymbol{T}_n\boldsymbol{B}_L)(i+1,j+1)=\boldsymbol{T}_n(i,j+1)+\boldsymbol{T}_n(i+1,j+2)+\boldsymbol{T}_n(i,j+2)$$
$$=t(i-1,j)+t(i,j+1)+t(i-1,j+1)=0$$

The sum is 0, for the three positions are just at the right place for the zero-sum rule Eq.(2), thus all the elements below the subdiagonal are 0s.

For the elements on the lower subdiagonal, $i-1=j\geq 0$, we have them all 1s :

$$(\boldsymbol{B}_L\boldsymbol{T}_n+\boldsymbol{T}_n\boldsymbol{B}_L+\boldsymbol{B}_L\boldsymbol{T}_n\boldsymbol{B}_L)(i+1,j+1)=t(i-1,j)+t(i,j+1)+0=b_j+b_{j+1}=1$$

In summary, the matrix product $\boldsymbol{S}_D\boldsymbol{S}_L$ is a bidiagonal matrix with both diagonals all 1s, or $\boldsymbol{S}_D\boldsymbol{S}_L=\boldsymbol{I}+\boldsymbol{B}_L$.

■

Such unipotent matrices can be used in matrix construction or factorization for Jordan canonical form blocks with prescribed eigenvalues or Frobenius normal form blocks with prescribed companion matrices, and the unipotent matrices are easy to generate recursively with integers only. By this theorem, we know that any unit triangular matrix is the product of two triangular unipotent matrices of index 2, any unipotent matrix can be a product of two unipotent matrices of index 2, and the factorization is not unique.

For example with Theorem 4, we use the idempotent matrix of $\boldsymbol{T}_4$ as we give in the proof of Theorem 1:

$$\boldsymbol{T}_4=\begin{bmatrix}a_0&0&0&0\\a_1&b_1&0&0\\a_2&-a_1-b_1&b_2&0\\a_3&a_1+b_1-a_2&a_1+b_1-b_2&b_3\end{bmatrix},\boldsymbol{B}_L=\begin{bmatrix}0&0&0&0\\1&0&0&0\\0&1&0&0\\0&0&1&0\end{bmatrix}$$

Then we have

$$\boldsymbol{N}_D=\boldsymbol{B}_L\boldsymbol{T}_n=\begin{bmatrix}0&0&0&0\\a_0&0&0&0\\a_1&b_1&0&0\\a_2&-a_1-b_1&b_2&0\end{bmatrix},\boldsymbol{N}_L=\boldsymbol{T}_n\boldsymbol{B}_L=\begin{bmatrix}0&0&0&0\\b_1&0&0&0\\-a_1-b_1&b_2&0&0\\a_1+b_1-a_2&a_1+b_1-b_2&b_3&0\end{bmatrix}$$

$$\boldsymbol{N}_D^2=\begin{bmatrix}0&0&0&0\\0&0&0&0\\a_0b_1&0&0&0\\-a_0(a_1+b_1)+a_1b_2&b_1b_2&0&0\end{bmatrix}=0,\boldsymbol{N}_L^2=\begin{bmatrix}0&0&0&0\\0&0&0&0\\b_1b_2&0&0&0\\b_1(a_1+b_1-b_2)-(a_1+b_1)b_3&b_2b_3&0&0\end{bmatrix}=0$$

$$\boldsymbol{S}_D\boldsymbol{S}_L=(\boldsymbol{I}+\boldsymbol{N}_D)(\boldsymbol{I}+\boldsymbol{N}_L)=\begin{bmatrix}1&0&0&0\\a_0&1&0&0\\a_1&b_1&1&0\\a_2&-a_1-b_1&b_2&1\end{bmatrix}\cdot\begin{bmatrix}1&0&0&0\\b_1&1&0&0\\-a_1-b_1&b_2&1&0\\a_1+b_1-a_2&a_1+b_1-b_2&b_3&1\end{bmatrix}$$

$$=\begin{bmatrix}1&0&0&0\\a_0+b_1&1&0&0\\(-1+b_1)b_1&b_1+b_2&1&0\\-(a_1+b_1)(-1+b_1+b_2)&(-1+b_2)b_2&b_2+b_3&1\end{bmatrix}=\begin{bmatrix}1&0&0&0\\1&1&0&0\\0&1&1&0\\0&0&1&1\end{bmatrix}=\boldsymbol{I}+\boldsymbol{B}_L$$

where we use $b_0=a_0, b_{i+1}=1-b_i, (1-b_i)b_i=0, b_i^2=b_i$.

## 9. Further Thoughts beyond

In the zero-sum triangles, in view of the general term of the interior cells Eq.(6), if the both boundary inputs are not all zeros, as long as they are nonzeros, positive or negative, or mixed with zeros, the further from the boundary edges, the bigger the cell numbers.

In the left Triangle 7 below, for the boundary input, only 0s on the left edge, $-1$s and 1s on the right edge, the total input is zero, but the further from the boundaries, in the middle of the 9th row, we have a cell with a number as big as 33.

```
                0
              0   -1
            0   1   1
          0   -1  -2  -1
        0   1   3   3   -1
      0   -1  -4  -6  -2   1
    0   1   5   10  8   1   1
  0   -1  -6 -15 -18  -9  -2  -1
0   1   7   21  33  27  11  3   1
```

Triangle 7

```
                 1
              -1   0
            -1   1   1
          -1   0  -2   0
        -1   1   2   2   1
      -1   0  -3  -4  -3   0
    -1   1   3   7   7   3   1
  -1   0  -4 -10 -14 -10  -4   0
-1   1   4  14  24  24  14   4   1
```

Triangle 8

In the right-hand Triangle 8, for the boundary input, only $-1$s on the left edge, 0s and 1s on the right edge, the total input is negative, but the further from the boundaries, in the middle of the 9th row, we have 2 cells with numbers as big as 24.

According to Stirling's approximation to the central and maximal binomial coefficient of the binomial distribution [16], the triangle numbers further from the boundaries can be increasing at a geometric rate.

If we go back to zero-sum's origin in economic game theory and view a zero-sum triangle through the lens of game theory, we can draw some interesting parallels between our matrix construction and a trading machine simulation. The boundary edge numbers can be seen as initial balances of “low-level” virtual traders following the zero-sum trading rule, while the positive and negative integers can be seen as the trader's cash surpluses and deficits, or credits and debts. Integers further away from the boundary edge can be seen “higher-level” traders' balances, which were determined by the trading patterns of “lower-level” traders. As we can see in Triangles 7&8, when trading moves further away from the starting positions at triangle boundaries, in terms of absolute magnitude, “higher-level” traders can see increasingly large balance surpluses and deficits, or like market bubbles of huge cash waves and financial crises of disastrous cash losses. This is the case even if the starting boundary inputs are small. This phenomenon lends itself to be compared with the real-world financial crises in 2008, where inter-linked financial institutions traded on highly-aggregated financial derivatives, whose default led to a disastrous outcome on the whole economy [17].

The zero-sum triangles may be a simple tool to construct integer balance triangles and to analyse the economic agents' trading behaviours that abide certain rules, so the potential areas of applications could include agent-based financial modelling and behaviour simulation of the economic agents in a zero-sum trading environment. The properties of the zero-sum triangles as triangular matrices should be able to find their explanations and applications in mathematical economics, which merits further investigation.

Every matrix property of the zero-sum triangles may lead to a number of combinatorial identities, compared to those related to Pascal's triangle [4], and even a specific triangle has some identity relations in combinatorics.

With the zero-sum rule, the proofs of the identities may also be much easier to clinch.

For the triangle construction itself, if three boundary edges are involved, a pyramid of cells can be built with the zero-sum rule, or more edges for higher dimensional pyramids. And the triangles can also be extended to rectangles by using negative indices and indices bigger than the boundary, and extending the triangle from the boundary edges outwards to exterior cells with the zero-sum rule, similar to that with Pascal's triangle [16].

**Acknowledgements**: This work was supported by the National Nature Science Foundation of China under Grant 62071013 and National Key R&D Program of China under Grant 2018AAA0100300.